\documentclass{article}
\usepackage{a4wide}
\usepackage{theorem}
\usepackage{amssymb}
\newcommand\al{\alpha}

\newcommand\ep{\varepsilon}

\newcommand\si{\sigma}

\newcommand\Ga{\Gamma}

\def\ZZ{\mathbb Z}

  \newcounter{numero}[section]
  \newcounter{aux1}
  \newcounter{aux2}
  \renewcommand{\thenumero}{\thesection.\arabic{numero}}
  \newenvironment{theorem}{
  \dimen255=\parindent \parindent=0in
  \refstepcounter{numero}{\vspace{.3cm}\bf Theorem \thenumero .}\
  \it}{\rm \parindent=\dimen255\par}

  \newenvironment{lemma}{
  \dimen255=\parindent \parindent=0in
  \refstepcounter{numero} {\vspace{.3cm}\bf Lemma \thenumero .}\
  \it}{\rm \parindent=\dimen255\par}

  \newenvironment{remark}{
  \dimen255=\parindent \parindent=0in
  \refstepcounter{numero} {\vspace{.3cm}\noindent\bf Remark \thenumero .}\
  }{\parindent=\dimen255
  \par}

  \newenvironment{proof}{
  {\vspace{.3cm}\noindent\sc{Proof:}}\
  }{$\fbox{}$ \par}

\usepackage{epsfig}

\title{Presentations for the monoids of 
singular braids \\ on closed surfaces}
\author{ Juan Gonz\'alez-Meneses}
\date{November, 2000}

\begin{document}
\maketitle

\begin{abstract}
 We give presentations, in terms of generators and relations, for the 
monoids $SB_n(M)$ of singular braids on closed surfaces. 
The proof of the validity of these presentations can 
also be applied to verify, in a new way, the presentations given 
by Birman for the monoids of Singular Artin braids.
\footnote{\hspace{-.65cm} Keywords: Braid - Singular Braid - Surface -  
Monoid - Presentation.}
\footnote{\hspace{-.65cm} {\em Mathematics Subject Classification:} 
 Primary: 20F36. Secondary: 20F05.}
\footnote{\hspace{-.65cm} Partially supported by 
DGESIC-PB97-0723 and by the european network TMR Sing. Eq. Diff. et Feuill.}
\end{abstract}

\section{Introduction}
 In this paper we deal with the braid groups of a closed surface
$M$. These groups are a natural generalization of Artin braid groups 
\cite{artin} and of the fundamental group of $M$. They are also subgroups 
of some {\em Mapping Class groups} of $M$, and finally they are fundamental
groups of the so called {\em Configuration spaces} of $M$ (see \cite{birman}
for a general exposition).

 They can be defined as follows. Fix $n$ ($n\geq 1$) distinct points 
$\{P_1,\ldots,P_n\}\in M$. A {\em n-braid} on $M$ is an n-tuple 
$b=(b_1,\ldots,b_n)$ of disjoint smooth paths $b_i$ in $M\times [0,1]$,
such that for all $i$, the path $b_i$ runs, monotonically on $t\in [0,1]$, 
from $(P_i,0)$ to some $(P_j,1)$. These $n$-braids are considered modulo
{\em isotopy} (deformation of braids fixing the ends), and there exists 
a multiplication of braids, given by concatenation of paths. The set of 
isotopy classes of $n$-braids on $M$, along with this multiplication,
forms the {\em braid group with $n$ strings on $M$}, denoted by $B_n(M)$.

 The following is a simple presentation of $B_n(M)$, in terms of generators and
relations, where $M$ is a closed, orientable surface of genus $g$ \cite{gm}:

\begin{itemize}
 
 \item Generators: \quad $\si_1, \ldots, \si_{n-1},  a_1, \ldots, a_{2g}. $
  
 \item Relations:
 
 \begin{itemize}
 
  \item[(R1)] $\si_i \si_j = \si_j \si_i$  \hspace{\stretch{1}}$(|i-j|\geq 2)$
 
  \item[(R2)] $\si_i \si_{i+1} \si_i = \si_{i+1} \si_i \si_{i+1}$
     \hspace{\stretch{1}}    $(1 \leq i \leq n-2)$
 
  \item[(R3)] $ a_1\cdots a_{2g}
        a_{1}^{-1}\cdots a_{2g}^{-1} =
        \si_{1}\cdots \si_{n-2} \si_{n-1}^2 \si_{n-2} \cdots \si_{1} $
 
  \item[(R4)] $  a_{r} A_{2,s}= A_{2,s} a_{r}$ \hspace{\stretch{1}}
    $(1\leq r,s \leq 2g; \;\;  r \neq s )$
 
  \item[(R5)] $ \left( a_{1}\cdots a_{r}\right) A_{2,r}=\si_1^2 A_{2,r}
      \left( a_{1}\cdots a_{r}\right)$  \hspace{\stretch{1}}
       $ (1\leq r\leq 2g) $
 
  \item[(R6)] $a_{r}\si_i=\si_i a_{r}$  \hspace{\stretch{1}}
   $ (1\leq r \leq 2g; \; \; i\geq 2)$
 
  \end{itemize}
 
\end{itemize}
where
$$
    A_{2,r}=\si_1^{-1}
   \left( a_{1}\cdots a_{r-1}a_{r+1}^{-1}\cdots a_{2g}^{-1}\right)
    \si_1^{-1}.
$$                  

\begin{figure}[ht]
\centerline{\input{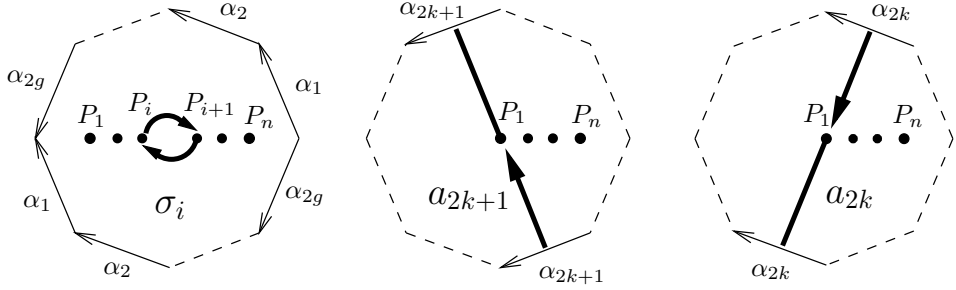}}
\caption{The generators of $B_n(M)$.}\label{generadores}
\end{figure}

The generators are represented in Figure~\ref{generadores}, where we have 
drawn the the canonical projections on $M$ of the considered braids, 
and $M$ is represented as a polygon of $4g$ sides, pairwise identified.

 We can also find in~\cite{gm} a similar presentation, when $M$ is a 
non-orientable, closed surface. 

 In the same way that singular Artin braids were defined (see~\cite{birman2})
to study Vassiliev invariants for these braids, we can define {\em singular 
braids} on $M$. Their definition is the same that the one of non-singular 
braids, but this time we allow a finite number of {\em singular points} 
(transverse intersection of two strings). The isotopy classes of
these singular braids, with the analogous multiplication, form the
{\em monoid of singular braids with $n$ strings on $M$}, denoted 
by $SB_n(M)$. This monoid is used in~\cite{gmparis} to define the Vassiliev 
invariants of braids on closed, orientable surfaces, proving, among other 
results, that these invariants classify these braids.

In~\cite{birman2} we can find presentations for $SB_n$, the monoids of
singular Artin braids, in terms of generators and relations. The main 
result of this paper is to give presentations for $SB_n(M)$. We will see,
as well, that the proof of this results furnishes a new proof of the
validity of the presentations in~\cite{birman2}.

\section{Presentation of $SB_n(M)$}

\begin{figure}[ht]
\centerline{\input{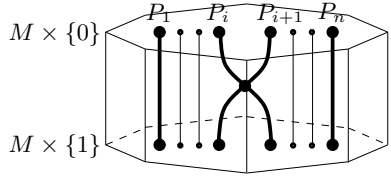}}
\caption{The singular braid $\tau_i$.}\label{taui}
\end{figure}

 We shall now give a presentation of $B_n(M)$, when $M$ is a closed, 
orientable surface of genus $g\geq 0$. The non-orientable case is
completely analogous, and is treated in a final remark at the end of this 
paper. We define, for all $i=1,\ldots,n-1$, the singular braid $\tau_i$ as 
in Figure~\ref{taui}, where the only non-trivial strings are the $i$-th
and the ($i+1$)-th ones, which intersect to form a singular point.
The result is the following:

\begin{theorem}\label{teosing}
 The monoid $SB_n(M)$ admits the following presentation:

\begin{itemize}
 
 \item Generators: \quad $\si_1, \ldots, \si_{n-1},  a_1, \ldots, a_{2g},
 \tau_1,\ldots, \tau_{n-1}$.
  
 \item Relations:
 
 \begin{itemize}
 
  \item[(R1-R6)] {\em Relations of $B_n(M)$}

  \item[(R7)] $\si_i \tau_j = \tau_j \si_i$
     \hspace{\stretch{1}}    $(|i-j| \geq 2)$
 
  \item[(R8)] $\tau_i \tau_j = \tau_j \tau_i$
     \hspace{\stretch{1}}    $(|i-j| \geq 2)$
 
  \item[(R9)] $\si_i \tau_i = \tau_i \si_i$
     \hspace{\stretch{1}}    $(i=1,\ldots, n-1)$
 
  \item[(R10)] $\si_i \si_j \tau_i = \tau_j \si_i \si_j$
     \hspace{\stretch{1}}    $(|i-j|=1)$
 
  \item[(R11)] $(a_{i,r} a_{i+1,r}) \tau_i (a_{i+1,r}^{-1} a_{i,r}^{-1})= 
\tau_i$   \hspace{\stretch{1}}    $(i=1,\ldots, n-1; \; r=1,\ldots,2g)$
 
  \item[(R12)] $\tau_i a_{j,r}= a_{j,r} \tau_i$
     \hspace{\stretch{1}}    $(j\neq i,i+1; \; r=1,\ldots,2g)$
 
  \end{itemize}
 
\end{itemize}
where
$$
a_{i,r}=\left\{\begin{array}{ll}
 (\si_{i-1}^{-1}\cdots  \si_{1}^{-1}) a_r (\si_{1}^{-1} \cdots \si_{i-1}^{-1})
& \mbox{ if $r$ is odd,} \\
 (\si_{i-1}\cdots  \si_{1}) a_r (\si_{1}\cdots \si_{i-1})
& \mbox{ if $r$ is even.}
\end{array} \right.
$$         
\end{theorem}

\begin{figure}[ht]
\centerline{\input{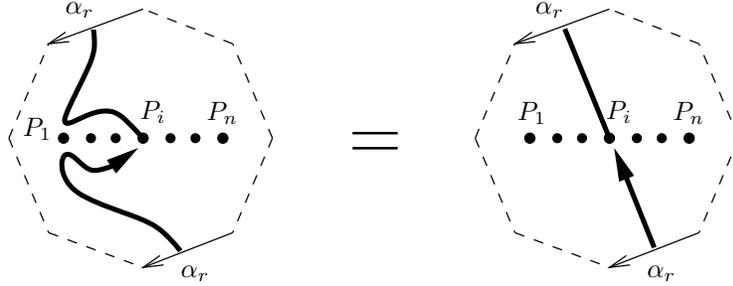}}
\caption{The braid $a_{i,r}=(\si_{i-1}^{-1}\cdots  \si_{1}^{-1}) a_r (\si_{1}^{-1} \cdots \si_{i-1}^{-1})$, when $r$ is odd.}\label{airair}
\end{figure}

\vspace{.3cm}
Remark that $a_{i,r}$ can be thought of as the $i$-th string 
crossing the ``wall'' $\al_r$, as we can see in Figure~\ref{airair} 
for the case when $r$ is odd.

\vspace{.3cm}
\noindent {\sc Proof of Theorem~\ref{teosing}:}
 First, it is evident that $\{\si_1, \ldots, \si_{n-1},  a_1, \ldots, a_{2g},
 \tau_1,\ldots, \tau_{n-1}\}$ is a set of generators of $SB_n(M)$, once
that we know (by \cite{gm}) that $\{\si_1, \ldots, \si_{n-1}, 
a_1, \ldots, a_{2g}\}$ generates $B_n(M)$.

\begin{figure}[ht]
\centerline{\input{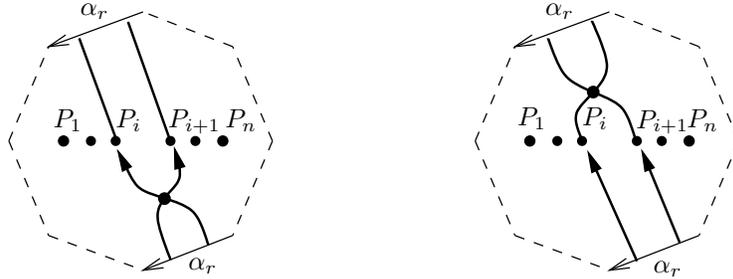}}
\caption{The braids $a_{i,r} a_{i+1,r} \tau_i$ and $\tau_i a_{i,r} a_{i+1,r}$
are isotopic, when $r$ is odd.}\label{R11}
\end{figure}

 It is also easy to prove that the proposed relations hold: (R1-R6) hold
in $B_n(M)$, which is a sub-monoid of $SB_n(M)$. (R7-R10) are known to hold
in $SB_n$, so they hold in a cylinder $D\times [0,1]$, where $D$ is a disk 
containing the $n$ points $P_1,\ldots,P_n$. We have just to extend the 
corresponding isotopy to all $M\times [0,1]$ by the identity.  
(R11) can be seen to hold in Figure~\ref{R11}, and finally (R12) is clear, since 
the only nontrivial strings of $\tau_i$ and $a_{j,r}$ can be isotoped to 
have disjoint projections on $M$, so these braids commute.

 In order to show that the relations are sufficient, we need the following 
lemma:

\begin{lemma}\label{cancela}
 The monoid $SB_n(M)$ is left-cancelative. That is, for all $a,b,c\in 
SB_n(M)$, one has: $c\:a=c\:b \Rightarrow a=b$.
\end{lemma}

\begin{proof}
 Since $\si_1, \ldots, \si_{n-1},  a_1, \ldots, a_{2g}$ are invertible, for
they belong to $B_n(M)$, we just need to prove that $\tau_i a= \tau_i b
\Rightarrow a=b$ for all $a,b \in SB_n(M)$ and all $i=1,\ldots,n-1$.

Thus, let us suppose that there exists an isotopy $H_t$ of $M\times[0,1]$,
such that $H_0=\mbox{id}_{M\times[0,1]}$ and $H_1(\tau_i a)=\tau_i b$.
Call $p$ the first singular point of $\tau_i a$ (the one corresponding
to $\tau_i$), and let $p_t=H_t(p)$. One has $p_0=p_1=p$.

 Let $V$ be the interior of a sphere of radius $\ep$ centered at $p$. We take 
$\ep$ small enough, such that $V\cap (\tau_i a)$ is as follows:

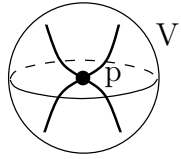
\begin{figure}[ht]
\centerline{\input{Vsing.pstex_t}}
\end{figure}

Denote $s_t=H_t(\tau_i a)$ and $V_t=H_t(V)$. We can suppose, without loss 
of generality, that $V_t$ is the interior of the sphere of radius $\ep$ 
centered at $p_t$, and that $V_t\cup s_t$ is as in the above picture.

Now, for $t\in[0,1]$, denote by $\widetilde{s}_t$ the braid which is obtained
by modifying $s_t$, only inside $V_t$, as follows:

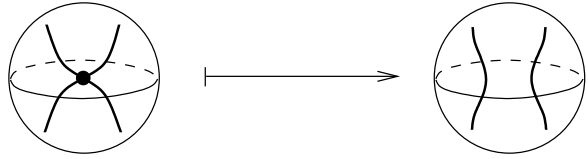
\begin{figure}[ht]
\centerline{\input{Vdesing.pstex_t}}
\end{figure}

We observe that $\widetilde{s}_0=a$ and $\widetilde{s}_1=b$, so
$H_t$ is an isotopy which transforms $a$ into $b$. Therefore, $a=b$. 
\end{proof}

\vspace{.3cm}
 Let us then show that Relations (R1-R12) are sufficient. 
Let $b,b'\in SB_n(M)$ be two isotopic singular braids, written in the 
generators of Theorem~\ref{teosing}. We must show that we can transform $b$
into $b'$ by using Relations (R1-R12). 

Being isotopic, both braids have the same number of singular points, say $k$. 
If $k=0$, the result follows from \cite{gm}, since (R1-R6) are sufficient 
relations for $B_n(M)$.

Suppose that $k>0$, and the result holds for braids with less than $k$
singular points. We can assume that the first letter of $b$ is 
$\tau_i$, for some $i$ (otherwise we can multiply $b$ and $b'$ on the left
by the greatest nonsingular ``prefix'' of $b$).
We will show that, using (R1-R12), we can transform $b'$ into a braid whose
first letter is $\tau_i$. The result then follows from Lemma~\ref{cancela},
and by induction hypothesis.

Let $p$ be the point of $b'$ corresponding (via isotopy) to the first 
singular point of $b$. This point $p$ must correspond to some $\tau_j$,
letter of $b'$. By (R10) and the braid relations (R1-R2), we can easily 
deduce the following:
$$
 \tau_j = \left\{ \begin{array}{ll}
 (\si_{j-1} \si_{j-2}\cdots \si_i)(\si_{j} \si_{j-1}\cdots \si_{i+1}) \tau_i
(\si_{i+1}^{-1}\cdots \si_{j}^{-1})(\si_{i}^{-1}\cdots \si_{j+1}^{-1})
& \mbox{ if } i<j,  \\
 (\si_{j+1} \si_{j+2}\cdots \si_i)(\si_{j} \si_{j+1}\cdots \si_{i-1}) \tau_i
(\si_{i-1}^{-1}\cdots \si_{j}^{-1})(\si_{i}^{-1}\cdots \si_{j+1}^{-1})
& \mbox{ if } i>j.  \\
\end{array}
\right.
$$
Hence, we can assume that the letter corresponding to $p$ is $\tau_i$.

Let us then write $b'=u\: \tau_i \: v$, where $u,v\in SB_n(M)$ and
$\tau_i$ is the above letter. Since $b$ is isotopic to $b'$, we can assume,
up to replacing $\tau_i$ by $\si_i \tau_i \si_i^{-1}$ (using (R9)), that
the $i$-th string of $u$ ends at the point $(P_i,s)$, for some $s\in[0,1]$. 
Hence, its canonical projection on $M$ is a loop in $M$ based at $P_i$, 
which induces an 
element $\mu\in\pi_1(M,P_1)$. This element can be modified as desired: 
it suffices to use (R11), replacing $\tau_i$ by 
$a_{i,r}^{\ep} a_{i+1,r}^{\ep} \tau_i \: a_{i+1,r}^{-\ep} a_{i,r}^{-\ep}$ 
($\ep=\pm 1$), to have $\mu$ transformed
into $\mu \:\overline{a}_{i,r}^{\ep}$, where $\overline{a}_{i,r}^{\ep}$ is
the projection on $M$ of the $i$-th string of $a_{i,r}^{\ep}$. Since
$\{\overline{a}_{i,1},\ldots, \overline{a}_{i,2g}\}$ is a set
of generators of $\pi_1(M,P_i)$, we can assume that $\mu=1$.

 Now notice that $u$ has less than $k$ singular points, hence any braid 
isotopic to $u$ can be obtained from it by applying (R1-R12), by 
induction hypothesis. We can then deform $u$ in such a way that its
$i$-th string will not go through the ``walls'' $\al_1, \ldots, \al_{2g}$
(we can do this since $\mu=1$).

 Let us go back to $b$, and consider a ``band'' $\Ga$, determined by
the $i$-th and the ($i+1$)-th strings of $b$, and which goes from $s=0$
to the first singular point of $b$, as in the figure below.

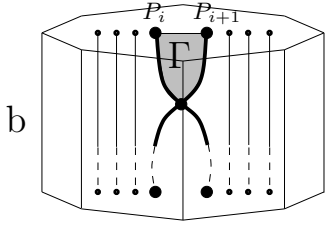
\begin{figure}[ht]
\centerline{\input{Gamma75.pstex_t}}
\end{figure}

Consider also an isotopy $H_t$ which transforms $b$ into $b'$. Recall
that $b'=u\: \tau_i v$, where the $i$-th string of $u$ does not go through
the walls. We can now consider $\Ga_1=H_1(\Ga)$, and deform the ($i+1$)-th 
string of $u$ along this band, in such a way that it will be as close to 
the $i$-th string as desired (recall that we are allowed to deform $u$).
We can then assume that neither the $i$-th nor the ($i+1$)-th string of 
$u$ goes through the ``walls'' of the cylinder $M\times [0,1]$.
Moreover, using (R9) we can modify the number of crossings of these two 
strings, as desired (just replacing $\tau_i$ by 
$\si_i^{r} \tau_i \si_i^{-r}$, $r\in\ZZ$). Therefore, we can assume that 
they do not cross, i.e. there is no $\si_j$ in $u$ involving the $i-th$
and the ($i+1$)-th strings. 

 We can also assume that these two strings are so close that one has the 
following property: if there is a letter $\si_j^{\ep}$
($\ep=\pm 1$) of $u$ which involves the $i$-th or the ($i+1$)-th string, 
then this letter, together with either the previous or the following one, 
forms a sub-word of $u$ of one of the following four types:

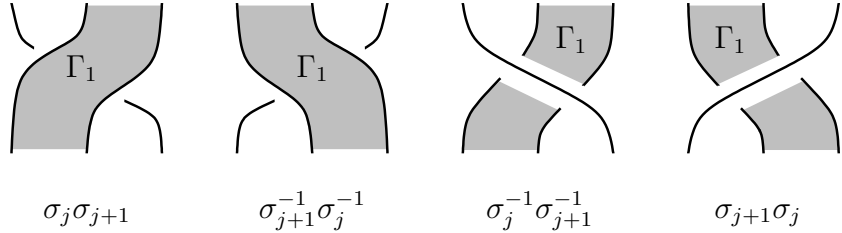
\begin{figure}[ht]
\centerline{\input{cruces.pstex_t}}
\end{figure}

 But in this case it is easy to see that, using relations (R7), (R8),
(R10) and (R12), we can ``raise'' the point $p$, until we get $\tau_i$ as
the first letter of $b'$. So by Lemma~\ref{cancela} we can cancel $\tau_i$,
and by induction hypothesis the resulting braids are equivalent by means
of Relations (R1-R12).
This ends the proof of Theorem~\ref{teosing}
$\fbox{}$

\begin{remark}
  There is an analogous presentation of $SB_n(M)$, when $M$ is a 
non-orientable, closed surface. We just need to consider the presentation 
given in~\cite{gm} for $B_n(M)$. Then replace, in the presentation of 
Theorem~\ref{teosing}, the generators $a_1,\ldots,a_{2g}$ by the 
corresponding generators on the non-orientable surface, and Relations
(R1-R6) by the relations given in~\cite{gm}. The same proof remains valid.
\end{remark}

\begin{remark}
 The presentation given in Theorem~\ref{teosing} can be easily simplified.
It suffices to eliminate the generators $\tau_2,\ldots,\tau_{n-1}$, replacing
in the relations $\tau_3$ by $(\si_2 \si_1 \si_3 \si_2) \tau_1 (\si_2^{-1}
 \si_3^{-1}\si_1^{-1} \si_2^{-1})$, and
eliminating all relations containing some $\tau_j$ ($j\neq 1,3$), since they
are obtained from the remaining ones. We proposed the presentation above
since it is more useful for handling singular braids.
\end{remark}

\begin{remark}
 We can also replace (R1-R6) by any other set of sufficient relations
for the given generators of $B_n(M)$.
\end{remark}

\begin{remark}
 The above proof of Theorem~\ref{teosing}, after eliminating every allusion to
$\pi_1(M)$, is a new proof of the validity of the presentation for
$SB_n$ proposed in~\cite{birman2}.

\end{remark}

\vspace{.3cm}
\noindent {\bf Acknowledgements:} I wish to thank Prof. Orlando Neto,
for giving me the opportunity to enjoy the excellent working conditions I 
found at the {\em Centro de Matem\'atica e Aplica\c coes Fundamentais} of the
University of Lisbon, where this paper was written down.

\vspace{.3cm}
\begin{flushleft}
\noindent J. GONZ\'ALEZ-MENESES \\
Departamento de Matem\'atica Aplicada I \\
Escuela T\'ecnica Superior de Arquitectura \\ 
Avda. Reina Mercedes, 2 \\ 
41012-Sevilla (Spain) \\
{\em meneses@cica.es}
\end{flushleft}




\end{document}